\documentclass[11pt,a4paper]{article}

\usepackage{amsmath}
\usepackage{amssymb}
\usepackage{amsfonts}
\usepackage{verbatim}
\usepackage{enumerate}
\usepackage{float}
\usepackage{hyperref}
\usepackage[T1]{fontenc}
\usepackage[utf8]{inputenc}
\usepackage{authblk}
\usepackage{a4}
\usepackage{float}
\usepackage{graphicx}
\usepackage{caption}
\usepackage{subcaption}

\newcommand{\fr}{\mbox{$f_r$}}
\newcommand{\thetar}{\phi}

\renewcommand{\Re}{\,\textrm{Re}} 
\renewcommand{\Im}{\,\textrm{Im}}

\title{A Visualizable, Constructive Proof of the Fundamental Theorem of Algebra, and a Parallel Polynomial Root Estimation Algorithm}
\date{\vspace{-5ex}}
\author[1]{\small Christopher Thron }
\author[2]{\small Jordan Barry }

\affil[1]{\small Texas A\& M University-Central Texas 79549 USA}
\affil[2]{\small Texas A\& M University-Central Texas 79549 USA}

\affil[1]{\small email address: thron@tamuct.edu}
\affil[2]{\small email address: jbarry@tamuct.edu}

\begin{document}

\maketitle

\noindent{\bf Abstract:}  
This paper presents an  alternative  proof of the Fundamental Theorem of Algebra that has several distinct advantages. The proof is based on simple ideas involving continuity and differentiation. Visual software demonstrations can be used to convey the gist of the proof. A rigorous version of the proof can be developed using only single-variable calculus and basic properties of complex numbers, but the technical details are somewhat involved.  In order to facilitate the reader's intuitive grasp of the proof, we first present the main points of the argument, which can be illustrated by  computer experiments. Next we fill in some of the details, using single-variable calculus. Finally, we give a numerical procedure for finding all roots of an $n$th degree polynomial  by solving $2n$ differential equations in parallel.\\

\noindent{\bf Keywords:} Fundamental Theorem of Algebra, calculus, chain rule, continuity.\\

\noindent{\bf AMS 2010 Subject Classification:} 30C15,65H04

\section{Introduction}\label{sec.Intro}
The fundamental theorem of algebra states that any polynomial function from the complex numbers to the complex numbers with complex coefficients has at least one root. There are several proofs of the fundamental theorem of algebra, which employ a number of different domains of mathematics, including complex analysis (Liouville's theorem, Cauchy's integral theorem or the mean value property) (\cite{SCHEP}\cite{VYBORNY}, topology (Brouwer's fixed point theorem)\cite{ARNOLD}, differential topology\cite{GUILLEMIN}, calculus (\cite{SCHEP},\cite{FEFFERMAN}), and ``elementary'' methods using meshes or lattices  \cite{ROSENBLOOM}, \cite{BRENNER}. For easily-accessible and readable web references that explain these proofs, see \cite{FILE},\cite{STEED},\cite{LINFORD},\cite{DUNFIELD}. Many of these proofs are beautiful and  elegant.  Most are not constructive and do not provide a practical method for  finding  roots (\cite{ROSENBLOOM} and \cite{BRENNER} are exceptions).  The proof we present here is both constructive, and provides a practical method for finding all roots of any polynomial through the numerical solution of differential equations with different initial conditions. Furthermore, we have created a simple, intuitive visual display that demonstrates the construction of roots.

\section{Empirical observations from computer experiments}

Consider the polynomial $f(z)=\sum_{n=0}^N a_n z^n$, where $N$ is a positive integer and $a_n$ are complex. We want to show that  $f(z)=0$ has at least one complex solution.  To approach this problem, we  make some preliminary empirical observations on the behavior of polynomial functions.

First we may consider how $f(z)$ behaves for some specific values of $z$. 
When $z=0$ we have $f(z)=a_0$, and when $z$  has a very large magnitude then the terms $a_n z^n$  in $f(z)$ also have large magnitudes, especially the leading-order term $a_N z^N$. To understand the behavior of $f(z)$ between these two extremes, we isolate the behavior of $f(z)$ for different values of $|z|$, as described below.  

A  complex number can be written in polar form as $z = r e^{i \theta}$, where $r>0$ is  the magnitude of $z$  and $0 \le \theta< 2\pi$. If we fix $r$ and allow $\theta$ to vary, then the set of points $\{ r e^{i\theta}, 0 \le \theta < 2\pi \}$ is a circle of radius $r$ in the complex plane, which we denote as $C_r$. Since the function $f$ is defined on all complex numbers, in particular it is defined on each circle $C_r$. The image of $C_r$ under $f$ is also a set (a curve, actually)  in the complex plane, which we may denote as $f(C_r)$.

\begin{figure}
\centering
  \includegraphics[width=6.0in]{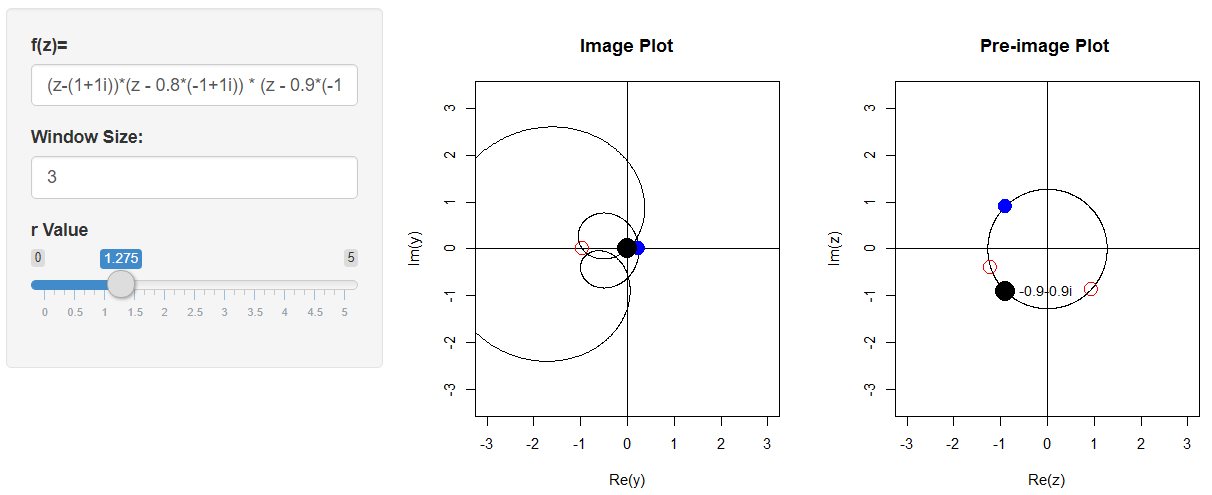}
  \caption{R Shiny interface for dynamic display of $f(C_r)$ and its preimage $C_r$, which shows upcrossings (solid blue dots) and downcrossings (hollow red dots). The function in this case is a cubic with roots at $1+i, -0.8+0.8i,$ and $-0.9 - 0.9i$. The large black dot in the preimage plot is the root $-0.9-0.9i$, which maps to 0 in the image plot. The R Shiny app used to create these plots is available online at \url{https://github.com/jthomasbarry/complex_plot_r}.}
  \label{fig:FTOA_shiny}
\end{figure}

\begin{figure}
\centering
  \includegraphics[width=5.0in]{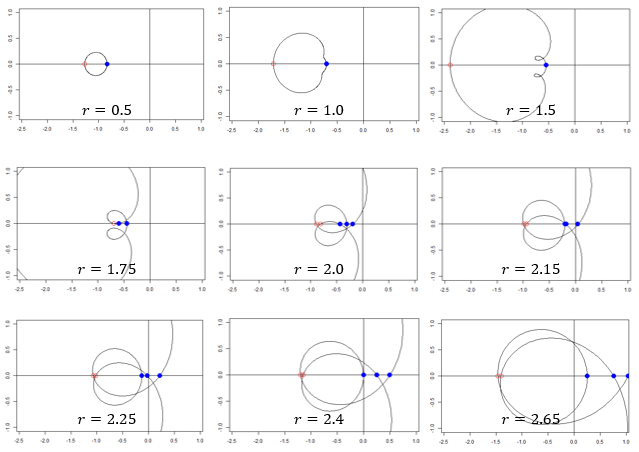}
  \caption{Curves $f(C_r)$ for different values of $r$, for the polynomial $f(z) =(a_1a_2a_3)^{-1}(z - a_1)(z-a_2)(z-a_3)$ where $a_1=1.6(1+i), a_2= 1.7(-1+i), a_3= 1.5(-1-i)$. The solid blue dots indicate upcrossings (which move to the right), while the hollow red dots indicate downcrossings (which move to the left). A new downcrossing-upcrossing pair is introduced  when $r\approx 1.75$ (as a lobe in the upper half plane expands down across the real axs) and another pair is introduced when $1.75 < r < 2$ (when a lobe in the lower half plane expands up across the real axis). Roots are found at upcrossings for moduli $r\approx 2.15, 2.25, 2.4$ (compare $|a_3|=2.12, |a_1| = 2.26, |a_2| = 2.40$).}
  \label{fig:FTOA_all}
\end{figure}

 Using computer software,  we may investigate the changes in the shape of $f(C_r)$ as $r$ increases from 0, for different polynomials $f(z)$.  
 For this purpose, an  R Shiny code (listed in the Appendix) was developed that  displays $f(C_r)$ and $C_r$ as well as  upcrossings and downcrossings,  for any given value of $r$  for an arbitrary  polynomial with complex coefficients, as specified by the user. A screenshot of the interface is shown in Figure~\ref{fig:FTOA_shiny}. A sequence of $f(C_r)$ plots for different values of $r$  is shown in Figure~\ref{fig:FTOA_all}.

  Without loss of generality we may assume $a_0=-1$: given any polynomial with $a_0 \neq 0$ we may obtain a polynomial with the same roots and having constant coefficient $-1$ by dividing by $-a_0$. If we look at several different polynomials $f$ and see how the curve $f(C_r)$ evolves as $r$ increases, we may make the following observations:
\begin{enumerate}[(i)]
\item
When $r$ is sufficiently small, then $f(C_r)$ has a nearly circular shape with center $-1$ and small  radius. The curve $f(C_r)$ has multiple intersections with the real axis.
\item
As $r$ increases, these points of intersection betwee $f(C_r)$ and the real axis move continuously along the real axis (although sometimes they disappear: see point (v) below)  
\item
There are two types of intersections: some move consistently to the right as $r$ increases, and others move consistently to the left.  When $r$ is small, the rightmost intersection is always rightward-moving.
\item
New  intersections with the $x$ axis may appear as $r$ increases. From a geometrical viewpoint, these new intersections occur when a lobe of  $f(C_r)$ located in the upper (resp. lower) half-plane shfts downward (resp. upward) as $r$ increases so that it intersects the axis.  These new intersections always appear first as a single point that splits into a left-moving and right-moving intersection as $r$ increases.  
\item
A right-moving intersection continues to move to the right unless it runs into a left-moving intersection, in which case both intersections may disappear. From the two-dimensional viewpoint, this occurs when a lobe of $f(C_r)$ that intersects the real axis moves above or below the axis.
\item
When $r$ is very large, the shape of $f(C_r)$ approaches a large circle centered at the origin. In particular, the rightmost intersection between $f(C_r)$ and the real axis is large and positive.
\item
The rightmost intersection when $r$ is small is always continuously connected to the rightmost intersection when $r$ is large by a series of right-moving or left-moving intersections.  Since the origin is on the real axis between these two intersections, the origin must be either a right-moving or left-moving intersection for some value of $r$.
\end{enumerate}

To understand the difference between right-moving and left-moving intersections, we may look more closely into the nature of the curve $f(C_r)$. As we mentioned above, the circle $C_r$ is parametrized by the angle $\theta$. As $\theta$ increases, the corresponding point on $C_r$ (given by $r e^{i \theta}$)  moves counterclockwise around $C_r$, while the image of the point under the function $f$ (given by $f(r e^{i\theta}$  traces out  the curve $f(C_r)$. As the tracing point crosses the real axis, we find there are two types of crossings: either upcrossings (from  below to above), or downcrossings (from above to below). It may be observed experimentally (and we shall soon show mathematically) that the upcrossings correspond to the rightward-moving  intersection points as noted above, and downcrossings correspond to leftward-moving intersection points.     
  
We may summarize a systematic procedure for using the software to locate  roots of $f$:  
\begin{enumerate}[(I)]
\item
Ensure that $a_0 = -1$ by dividing $f$ by $-a_0$  (if $a_0=0$, then 0 is  a root already);
\item  
Start with a small value of $r$ and locate an upcrossing point on $C_r$;  
\item
Follow the upcrossing point as it moves rightward. Eventually it will either pass over the origin, or run into a leftward-moving downcrossing point and disappear.
\item 
If the latter holds, follow  the leftward moving point backwards (i.e. decreasing $r$). Eventually, either it will pass over the origin, or it will merge with a rightward-moving point and disappear.
\item
Follow this rightward-moving point forward (increasing $r$) until it either passes through the origin or merges with a leftward-moving point.
\item
Continue iterating Steps IV and V  until the origin is reached. 
\end{enumerate}

\section{Outline of a formal proof}

This procedure is the basis for a formal proof of the theorem.  Some of the technical detals are rather involved, but the guiding intuition  is captured by the procedure described above. 

The proof proceeds in several steps:

\begin{enumerate}
\item
The roots of $f(z)$ are  identical to the roots of $- f(z)/a_0$.  So without loss of generality, we may assume that the constant coefficient $a_0$  is equal to $-1$.
\item
Assume for the moment that $f^{\prime}(z)$ has no zeros on the real axis between $-1$ and $0$.  (Later we will deal with the case where this is not true.) 
\item
For any value of $r>0$, we define the curve $\fr (t) \equiv f(re^{it}), 0 \le  t \le  2\pi$.  Since   $\fr (0)=\fr (2\pi)$, it follows that this is a closed (possibly self-intersecting) curve in the complex plane.
\item
Denote by  \emph{upcrossing} (resp. \emph{downcrossing}) a point where $\fr$ crosses the real axis from below (resp. above).  In other words,  the real number $x$ is an upcrossing for $\fr$ if there exists $t$ such that $\fr (t)=x$, and there exists $\delta>0$  such that  $\Im \fr (s) \le 0$ for $t - \delta \le s \le t$  and $\Im \fr (s) \ge 0$ for $t \le s \le t + \delta$. By continuity, every root of $\Im \fr$ is either an upcrossing, a downcrossing, or a point where the real line is tangent to $\fr$.
\item
Suppose $x_r = \fr(t)$ is an upcrossing and $\fr^{\prime}(t) \neq 0$, then the crossing point moves continuously to the right  as a function of $r$. More precisely, there exist $\epsilon, \delta > 0$   and a continuous, real-valued function $g(s)$ that is strictly increasing on   the interval $r - \epsilon < s <  r + \epsilon$  such that $g(r) = x_r$ and $g(s) = f_s(t')$, for some $t'$  in the interval $t-\delta < t' < t+\delta$.  We call the function $g$ an \emph{upcrossing function}.
A similar statement holds for the downcrossing case, except that $g$ is decreasing: the function $g$ in this case is called a \emph{downcrossing function}).
\item
The domain of any upcrossing function $g$ may be extended to an open interval, such that  either the range of $g$ includes the origin, or the right endpoint $b$ of the domain is such that $g(b)$  is a point of tangency of the curve $f_b$.   
\item
Every point of tangency that is the right endpoint of the domain of an upcrossing function  is the left endpoint of  the domain of a downcrossing function.
\item 
Every point in the interval $[-1,0]$ is either an upcrossing, a downcrossing, or a point of tangency of $\fr$ for some positive value of $r$.  In particular, the origin is either an upcrossing, downcrossing, or point of tangency, and is thus equal to $f(r e^{it})$ for some values of $r$ and $t$.
\end{enumerate}

Steps (1-8) handle the case where none of the roots  of $f^{\prime}$ lie on the real segment $[-1,0]$. If on the other hand $f^{\prime}$ does have a root on $[-1,0]$, we may consider the ray $\theta  = \pi + \nu$, for sufficiently small $\nu$, which (by continuity)  will have at least one upcrossing intersection with $C_r$   if $r,\nu$ are sufficiently small. We denote this intersection as $\widetilde{z}$. Since the roots of $f^{\prime}$ are isolated, we can also choose the $\nu$  such that $f^{\prime}$ has no roots on the ray. We may then consider the function
$\widetilde{f}(z) \equiv f(z) e^{-i \nu}$.
Then $\widetilde{z}e^{-i \nu}$ is an upcrossing point for $\widetilde{f}$ of the negative real axis. Steps (1-8) above then goes through for $\widetilde{f}$: and the roots of $\widetilde{f}$ are identical with the roots of $f$.  
    
\section{Parallel numerical  procedure for finding all roots of a polynomial}
As above, we suppose $f(z)= \sum_{n=1}^N a_n z^n$ with $a_0=-1$, where $z = r e^{i\theta}$.
We seek equations satisfied by upcrossing locations $x$ as a function of $r$. Note that in order for $x=f(re^{i\theta})$ to be an upcrossing, the complex argument $\theta$ varies as $r$ varies, so we must consider
both $\theta$ and $x$ as functions of $r$. To make this clear, we will use $\thetar = \thetar(r)$ to denote the complex argument, so that $x(r) = f(z(r))$ where $z(r)=r e^{i\thetar}$.

Using the chain rule, we have:
\begin{equation}\label{eq:dx_dr}
\frac{dx}{dr}= f^{\prime} \left(z\right) \frac{d}{dr} \left(z\right) = f^{\prime}\left(z\right) e^{i\thetar}\left(1 + ir\frac{d\thetar}{dr}\right)
\end{equation}
Since $x$ is real-valued, so $\frac{dx}{dr}$ is also a real function and  $\frac{dx}{dr} = \left(\frac{dx}{dr}\right)^*$, where $^*$ denotes complex conjugate. This gives
\begin{equation}
f^{\prime}\left(z\right) e^{i\thetar}\left(1 + ir\frac{d\thetar}{dr}\right) = f^{\prime}\left(z\right)^* e^{-i\thetar}\left(1 - ir\frac{d\thetar}{dr}\right).
\end{equation}
Solving for $r \frac{d\thetar}{dr}$, we obtain
\begin{equation}\label{eq:dth_dr}
r \frac{d\thetar}{dr}= \frac{ \left(f^{\prime}\left(z\right)^* e^{-i\thetar} -f^{\prime} \left(z\right) e^{i\thetar} \right)  } { i(f^{\prime} \left(z\right) e^{i\thetar} + f^{\prime}\left(z\right)^* e^{-i\thetar})}
= \frac{-\Im \left(f^{\prime} \left(r e^{i\thetar(r)}\right) e^{i\thetar(r)}\right)} 
{\Re \left(f^{\prime} \left(r e^{i\thetar(r)}\right) e^{i\thetar(r)}\right) },  
\end{equation}
where we have replaced $z$ with $r e^{i\thetar(r)}$ to highlight the $r$ dependence.
From \eqref{eq:dx_dr} and \eqref{eq:dth_dr} we may calculate:
\begin{equation}\label{eq:dx_dr2}  
\frac{dx}{dr}=  f^{\prime}\left(z\right) e^{i\thetar}\left(1 - i\frac{\Im(f^{\prime} \left(z\right) e^{i\thetar}} {\Re(f^{\prime} \left(z\right) e^{i\thetar} }\right)
=\frac{|f^{\prime} (re^{i\thetar(r)})e^{i\thetar(r)}|^2}
{\Re \left(f^{\prime} \left(r e^{i\thetar(r)}\right) e^{i\thetar(r)}\right) }  
\end{equation} 
Equations \eqref{eq:dth_dr} and \eqref{eq:dx_dr2} express $\frac{d\thetar}{dr}$ and $\frac{dx}{dr}$ respectively in terms of  $f^{\prime} \left(r e^{i\thetar(r)}\right) e^{i\thetar(r)}$. Fortunately, this rather complicated expression turns out to have a relatively simple interpretation.
The definition of $\fr$ implies that $\frac{d}{d\theta} \fr(\theta) = f^{\prime}(r e^{i\theta}) (i re^{i\theta})$, so if we define:
\begin{equation}\label{eq:def_alpha}
\alpha(r) \equiv \left. \frac{d}{d\theta} \fr(\theta)  \right|_{\theta=\thetar(r)}
\end{equation}
then we may re-express the system \eqref{eq:dth_dr}-\eqref{eq:dx_dr2} as:
\begin{equation}\label{eq:r_sys}
\begin{aligned}
 \frac{d\thetar}{dr} =\frac{-\Im \left(-i\alpha(r) \right)} { r \Re \left(-i\alpha(r) \right)}& =\frac{\Re \left(\alpha(r) \right)} { r \Im \left(\alpha(r) \right)};\\
\frac{dx}{dr} =\frac{|-i\alpha(r)|^2}{ \Re \left(-i\alpha(r) \right) }&=\frac{|\alpha(r)|^2}{ r \Im \left(\alpha(r) \right) }.
\end{aligned}
\end{equation}
For future reference, note that $\Im(\alpha(r))$ is positive or negative depending on whether $x(r)$ is an upcrossing or downcrossing.

Alternatively, we can pose the system such that $\theta$ is the independent variable, and $r,x$ are the dependent variables. To clarify the dependence of $r$ on 
$\theta$, we use $\rho = \rho(\theta)$ here to represent the complex modulus as a function of $\theta$, so that $f_{\rho}(\rho e^{i\theta})$ is an upcrossing point for the curve $C_{\rho}$. In analogy to \eqref{eq:def_alpha}, we define:
\begin{equation}\label{eq:def_beta}
\beta(\theta) \equiv  \left. \frac{d}{d\theta} f_{r} (\theta) \right|_{r=\rho(\theta)},
\end{equation}
and in analogy to \eqref{eq:r_sys} we obtain:
\begin{equation}\label{eq:theta_sys}
\begin{aligned}
 \frac{d\rho}{d\theta} =\frac{-r \, \Re \left(-i\beta(\theta) \right)} { r\Im \left(-i\beta(\theta) \right) } 
 &=\frac{\Im \left(\beta(\theta) \right)} { \Re \left(\beta(\theta) \right) };\\
\frac{dx}{d \theta} =-\frac{|-i\beta(\theta)|^2}{r \Im \left(-i\beta(\theta)\right) }
& =\frac{|\beta(\theta)|^2}{ r \, \Re \left(\beta(\theta)\right) }.
\end{aligned}
\end{equation}
It follows that given a crossing point $x = f_r(\theta)$, we can always make the crossing point `move' continuously to the right by following this strategy:
\begin{enumerate}
\item
If $\lim_{z \rightarrow r e^{i\theta}} | \Im f'(re^{i\theta})/\Re f'(re^{i\theta}) |> c$  (where $c < 1$ is a fixed positive parameter), then propagate $x$  to the right using \eqref{eq:r_sys}
 with either increasing or decreasing $r$, depending on the sign of $\Im f'(re^{i\theta})$:
 \item
 Otherwise, propagate $x$ to the right using \eqref{eq:theta_sys} with either increasing or decreasing $\theta$, depending on the sign of  $\Im f'(re^{i\theta})$.
\end{enumerate}
Following this procedure will yield a monotonically increasing  crossing point.  If the initial crossing point is chosen such that it is chosen between $-1$ and 0 on the real axis, then eventually the crossing point will pass 0 and a root will be obtained. 

The above procedure is only guaranteed to obtain a single root $\gamma_1$. Subsequent roots may be estimated by taking $f^{(1)}(z) \equiv f(z)(1-z/\gamma_1)^{-1}$ and finding another root $\gamma_2$, then iterating the procedure with $f^{(j)}(z) \equiv f(z)(1-z/\gamma_j)^{-1}, j = 1,2,\ldots$ until all roots are found.  However, it is possible there may be numerical stability problems, because due to numerical error $f^{(j)}(z)$ is no longer a polynomial for $j \ge 1$. 

An alternative approach finds all roots in parallel as follows. If the polynomial has degree $n$, the $z^n$ term dominates the behavior of $f(C_r)$ when $r$ is large. It follows that for $r$ sufficiently large,  $f(C_r)$ must have at least $n$ upcrossings on the positive real axis and at least $n$ downcrossings of the negative real axis.  All upcrossings may be followed leftwards using the reverse of the rightward-tracking procedure described above; and all downcrossings may be followed rightward by a similar procedure. Not all of these $2n$ tracks (which may be computed in parallel) will result in a root; however, it is guaranteed that all $n$ roots will be obtained through the procedure.


\begin{thebibliography}{999}
\bibitem{SCHEP}
Schep, A.: A simple complex analysis and and advanced calculus proof of the fundamental theorem of algebra. American Mathematical Monthly 116(1), 67--68 (2009).
\bibitem{VYBORNY}
Vyborny, R.: A simple proof of the fundamental theorem of algebra. Mathematica Bohemica 135(1), 57--61 (2010).
\bibitem{ARNOLD}
Arnold, B.H: A topological proof of the fundamental theorem of algebra. The American Mathematical Monthly 56(7), 465--466 (1949).
\bibitem{GUILLEMIN}
Guillemin, V., Pollack, A.: Differential Topology. 1st edn. Prentice-Hall (1947).
\bibitem{FEFFERMAN}
Fefferman, C.: An easy proof of the fundamental theorem of algebra. The American Mathematical Monthly 74(7),854--855 (1967).
\bibitem{ROSENBLOOM}
Rosenbloom, P.C.: An elementary constructive proof of the fundamental theorem of algebra.The American Mathematical Monthly 52(10), 562--570 (1945).
\bibitem{BRENNER}
Brenner, J.L, Lyndon, R.C.: Proof of the fundamental theorem of algebra. American Mathematical Monthly 88(4), 253--256 (1981).
\bibitem{FILE}
File, D., Miller, S.: Fundamental theorem of algebra lecture notes from the reading classics (euler) working group autumn 2003, \url{https://people.math.osu.edu/sinnott.1/ReadingClassics/FundThmAlg_DFile.pdf}, last accessed 2020/2/1.
\bibitem{STEED}
Steed, M.: Proofs of the fundamental theorem of algebra, \url{http://math.uchicago.edu/~may/REU2014/REUPapers/Steed.pdf}, last accessed 2020/1/12.
\bibitem{LINFORD}
Linford, K.: An analysis of Charles Fefferman's proof of the fundamental theorem of algebra, http://commons.emich.edu/honors/504, last accessed 2020/2/1.
\bibitem{DUNFIELD}
Dunfield, N.: The fundamental theorem of algebra (class notes), \url{https://faculty.math.illinois.edu/~nmd/classes/2015/418/notes/fund_thm_alg.pdf}, last accessed 2020/1/15.
\end{thebibliography}
\end{document}